\documentclass[12pt,leqno]{article}

\newtheorem{theorem}{Theorem}
\newtheorem{proposition}{Proposition}
\newtheorem{corollary}{Corollary}
\newtheorem{remark}{Remark}
\newtheorem{lemma}{Lemma}
\newfont{\bb}{msbm10 at 12pt}

\def\C{\hbox{\bb C}}
\def\R{\hbox{\bb R}}
\def\S{\hbox{\bb S}}
\def\P{\hbox{\bb P}}
\def\H{\hbox{\bb H}}

\begin{document}

\title{A Willmore functional for compact surfaces of complex projective plane}       

\author{Sebasti\'{a}n Montiel \thanks{Research partially 
supported by a DGICYT grant No. PB97-0785.} \and Francisco Urbano$^{\ *}$}

\date{February 2000}        

\maketitle
 
\section{Introduction}
Amongst the global conformal invariants for compact surfaces in a Riemannian manifold $(M,\langle,\rangle)$, perhaps the better known is the Willmore functional. For an immersion $\phi:\Sigma\rightarrow M$ from a compact surface $\Sigma$, the Willmore functional is defined by
\[
W(\phi)=\int_{\Sigma}\left(|H|^2+\bar{K}\right)dA,
\]
where $H$ denotes the mean curvature vector of the immersion $\phi$, $\bar{K}$ the sectional curvature of $M$ restricted to $\Sigma$ and $dA$ the canonical measure of the induced metric.

This functional has been extensively studied when $M$ is the Euclidean space $\R^n$ (or the sphere $\S^n$ or the hyperbolic space $\R\H^n$, because $W$ is invariant under conformal transformations of ambient space). The type of problems which have been studied for this ambient space $\R^n$ are of a different nature. On the one hand are the papers written with the object of obtaining lower bounds for the functional $W$. Given that $W$ has $4\pi$ as an absolute minimum and is only reached by the round spheres of $\R^n$, Willmore proposed to study this functional on tori and conjectured that $2\pi^2$ is the minimum value for these surfaces and it is attained only by the Clifford torus. This problem, which still has not been resolved, can be considered as the starting point of a series of important papers in which minimization problems of $W$ are studied. Amongst those which we mention are [K1], [K2], [K3], [LY], [M], [MoR], [R] and [S].

On the other hand, other authors are interested in the study of critical surfaces for $W$ (the known {\em Willmore surfaces}) whose Euler-Lagrange equations were obtained (only when $M$ has constant curvature) by Weiner in [W1]. Two papers of interest are [BB] and [P], where the authors construct Willmore tori in $\S^3$ from minimal surfaces of $\R\H^3$ and elastic curves of $\S^2$ respectively. Two other papers of interest are [B] and [Mo] where Willmore spheres in $\S^3$ and $\S^4$ respectively are classified, describing the relation of these Willmore spheres and certain class of minimal surfaces in $\R^3$ and $\R^4$. 

In this paper the authors studied the Willmore functional for compact surfaces of the complex projective plane $\C\P^2$. If $\phi:\Sigma\rightarrow\C\P^2$ is an immersion of an orientable compact surface $\Sigma$ in the complex projective plane of constant holomorphic sectional curvature $4$, then the Willmore functional is given by
\[
W(\phi)=\int_{\Sigma}\left(|H|^2+1+3C^2\right)dA,
\]
where $C$ is the K\"{a}hler function on $\Sigma$ defined by $\phi^*(\Omega)=C\,dA$, $\Omega$ being the K\"{a}hler two form on $\C\P^2$ and $dA$ the volume two form on $\Sigma$.

This functional can be written (see section 2) as $W={1\over 2}(W^++W^-)$, where $W^+$ and $W^-$ are also conformal invariant functionals (see Proposition 1), defined by
\[
W^+(\phi)=\int_{\Sigma}\left(|H|^2+6C^2\right)dA,\quad   W^-(\phi)=\int_{\Sigma}\left(|H|^2+2\right)dA.
\]
These functionals $W^{\pm}$ are closely related with the Penrose twistor bundles ${\cal P}^{\pm}$ over $\C\P^2$ (see Proposition 3), because twistor holomorphic surfaces, i.e., surfaces of $\C\P^2$ whose twistor liftings are holomorphic (which we study in depth in Theorem 1) are critical surfaces for these functionals (in fact they are minimizers for $W^{\pm}$). As $\C\P^2$ with its canonical orientation is a self-dual Riemannian manifold but not an anti-self-dual, the twistor bundle ${\cal P}^-$ is a complex manifold and the twistor bundle ${\cal P}^+$ is an almost-complex manifold but non-complex (see [AHS]). This fact allows us to easily construct (see Proposition 6) twistor holomorphic compact surfaces with negative spin (and therefore Willmore surfaces for $W^-$), but it is really complicated to get non-trivial examples of twistor holomorphic surfaces with positive spin. 

The Euler-Lagrange equations for the functionals $W^{\pm}$ (Proposition 5) say that the minimal surfaces of $\C\P^2$ are critical for the functional $W^-$. It is also interesting to remark that the functional $W^-$ restricted to minimal surfaces is twice the area functional. Due to all these considerations about $W^{\pm}$, the authors think that the $W^-$ Willmore functional is the natural one to be studied for surfaces in $\C\P^2$. In this way, in section 4, we study lower bounds for $W^-$ obtaining results which can be sumarized as follows

\begin{quote}

{\em Let $\phi:\Sigma\rightarrow\C\P^2$ be an immersion of a compact surface $\Sigma$. Then
\begin{description}
\item [{i)}] $W^-(\phi)\geq 2\pi\mu$, being $\mu$ the maximum multiplicity of $\phi$, and the equality holds if and only if $\phi(\Sigma)$ is a complex projective line. In particular $2\pi$ is the minimum value for $W^-$.
\item [{ii)}] If $\phi$ is Lagrangian, then $W^-(\phi)\geq 4\pi\mu$, and the equality holds if and only if $\phi(\Sigma)$ is either a real projective line or a Whitney sphere. In this case $4\pi$ is the minimum value for $W^-$.
\end{description}}
\end{quote}

Finally, we wish to mention that the proofs of the above results (Theorems $2$ and $3$ in the paper) can be easily extended for surfaces in the complex projective space $\C\P^n$, so that Theorems $2$ and $3$ and Corollaries $1$ and $2$ are also true when we change $\C\P^2$ by $\C\P^n$. We also wish to mention that almost all the results in the paper can be generalized for surfaces of the complex hyperbolic plane.

\section{The Willmore functional for surfaces in four manifolds}
Let $(M,\langle ,\rangle )$ be an n-dimensional Riemannian manifold and $\phi:\Sigma\rightarrow M$ an immersion of a compact surface $\Sigma$. The Willmore functional $W(\phi)$ is defined by
\[
W(\phi)=\int_{\Sigma}\left(|H|^2+\bar K\right)dA,
\]
where $H$ is the mean curvature of $\phi$, $dA$ the area two form on $\Sigma$ and $\bar K=\bar R (e_1,e_2,e_2,e_1)$, being $\bar R$ the curvature of $\langle ,\rangle $ and $\{e_1,e_2\}$ an orthonormal basis in $\Sigma$. This functional is invariant under conformal changes of the metric $\langle,\rangle$.

From now on we suppose that $M$ is an oriented four--manifold, and that $\Sigma$ is an oriented surface. If $T^{\perp}\Sigma$ is the normal bundle of $\phi$, then we have the orthogonal decomposition
\[
\phi^*TM=T\Sigma\oplus T^{\perp}\Sigma.
\]
Let $\bar {\nabla}$ be the connection on $\phi^*TM$ induced by the Levi-Civita connection of $TM$ and let $\bar{\nabla}=\nabla +\nabla^{\perp}$ be the corresponding decomposition.

If $\{e_1,e_2,e_3,e_4\}$ is an oriented orthonormal local reference on $\phi^*TM$ such that $\{e_1,e_2\}$ is an oriented reference on $T\Sigma$, then we define the {\em normal curvature} $K^{\perp}$ of the immersion $\phi$ by
\[ 
K^{\perp}=R^{\perp}(e_1,e_2,e_3,e_4),
\]
where $R^{\perp}$ is the curvature tensor of the normal connection $\nabla^{\perp}$. 
Also we will denote by $\bar K^{\perp}$ the function on $\Sigma$ given by
\[
\bar K^{\perp}=\bar R(e_1,e_2,e_3,e_4).
\] 
When $\Sigma$ is compact, the Euler characteristics of $T\Sigma$ and $T^{\perp}\Sigma$ are given respectively by
\[
\chi={1\over 2\pi}\int_{\Sigma}KdA  \quad {\rm and} \quad  \chi^{\perp}={1\over 2\pi}\int_{\Sigma}K^{\perp}dA.
\]
Now the Willmore functional can be decomposed as
\[
W(\phi)={1\over 2}\left(W^+(\phi)+W^-(\phi)\right),
\]
where 
\begin{eqnarray*}
&{\displaystyle W^+(\phi)=\int_{\Sigma}\left(|H|^2+\bar K -\bar K^{\perp}\right)dA,}& \\
&{\displaystyle W^-(\phi)=\int_{\Sigma}\left(|H|^2+\bar K+\bar K^{\perp}\right)dA.}&
\end{eqnarray*}

\begin{proposition}
The functionals $W^+$ and $W^-$ are invariant under conformal changes of the metric $\langle ,\rangle $ on $M$.
\end{proposition}
{\em Proof:\/}
Let $\langle, \rangle _*=e^{2u}\langle,\rangle $ a metric on $M$ conformal to $\langle,\rangle $, being $u:M\rightarrow \R$ a smooth function. Then it is very well-known that the second fundamental forms $\sigma$ and $\sigma_*$ of $\phi$ with respect to $\langle,\rangle$  and $\langle,\rangle _*$ are related by
\[
\sigma_*(u,v)=\sigma(u,v)-\langle u,v\rangle (\bar{\nabla}u)^{\perp},
\]
where $^{\perp}$ means normal component.
From here, it is an easy exercise to check that
\begin{eqnarray*}
&\left(|H_*|_*^2+\bar K_*-K_*\right)dA_*=\left(|H|^2+\bar K -K\right)dA,& \\
&\left(\bar K^{\perp}_*-K^{\perp}_*\right)dA_*=\left(\bar K^{\perp}-K^{\perp}\right)dA&
\end{eqnarray*}
where $*$ means the corresponding object for the metric $\langle,\rangle _*$. The above formulae prove the Proposition taking into account the Gauss-Bonnet theorem and that the normal bundles of $\phi$ (with respect to both metrics) are isomorphic.\hfill q.e.d.

Now we are going to relate these functionals with the twistor bundles. Given a point $x\in M$, let ${\cal P}_x^{\pm}$ be the set of almost Hermitian structures $J_x^{\pm}$ over $T_xM$ such that if $\Omega^{\pm}(u,v)=\langle J_x^{\pm}u,v\rangle $, then ${\pm}\Omega\wedge\Omega$ is the orientation induced on $T_xM$ from $M$. Then ${\cal P}^{\pm}=\cup_{x\in M}{\cal P}_x^{\pm}$ are $\C\P^1$-fiber bundles over $M$, called the {\bf twistor bundles} of $M$. If $\pi^{\pm}:P^{\pm}\rightarrow M$ are the projections, then the vertical distributions $V^{\pm}=\ker \pi^{\pm}_*$ inherit from the standard complex structure of the complex projective line $\C\P^1$ an almost complex structure $J^{v\pm}$. The Levi-Civita connection of $M$ induces a decomposition of the tangent bundles of ${\cal P}^{\pm}$
\[
T{\cal P}^{\pm}=H^{\pm}\oplus V^{\pm}
\]
with $H^{\pm}\equiv TM$, via $\pi^\pm_*$ and there are also on the horizontal distributions $H^{\pm}$  almost complex structures $J^{h\pm}$ defined by
\[
J^{h\pm}_{J^{\pm}_x}=J^{\pm}_x.
\]
So ${\cal J}^{\pm}=J^{h\pm}+J^{v\pm}$ define almost complex structures on ${\cal P}^{\pm}$ which only depend of the oriented conformal structure of $M$. A central result due to Atiyah, Hitchin and Singer [AHS] is that {\em $({\cal P}^+,{\cal J}^+)$ is a complex manifold if and only if $(M,\langle,\rangle )$ is anti-self-dual, and $({\cal P}^-,{\cal J}^-)$ is a complex manifold if and only if $(M,\langle,\rangle )$ is self-dual}.

Let $\phi:\Sigma\rightarrow M$ an immersion of an {\em oriented} surface $\Sigma$ and $\{e_1,e_2,e_3,e_4\}$ an orthonormal local reference on $\phi^*TM$ such that $\{e_1,e_2\}$ is an oriented reference on $T\Sigma$. We are going to define two almost complex structures $J^{\pm}$ on $\phi^*TM$ by
\[
J^{\pm}(e_1)=e_2, J^{\pm}(e_3)={\pm}e_4.
\]
We remark that the $J^{\pm}$ on $\Sigma$ is the complex structure on $\Sigma$ compatible with the given orientation, and that on $T^{\perp}\Sigma$, $\nabla^{\perp}J^{\pm}=0$. So the Koszul-Malgrange theorem [KM] says that $J^{\pm}$ give to $T^{\perp}\Sigma$ two unique structures of holomorphic line bundles over $\Sigma$ such that a normal section $\xi$ is holomorphic if and only if
\[
\nabla^{\perp}_{J^{\pm}v}\xi=J^{\pm}\nabla_v\xi
\]
for any $v\in T\Sigma$.

We define the {\em twistor liftings} $\tilde{\phi}^{\pm}:\Sigma\rightarrow {\cal P}^{\pm}$ by
\[
\tilde{\phi}^{\pm}(p)=J^{\pm}_{\phi(p)}.
\]
Although it is not explicity stated, the following result was proved in [F].

\begin{proposition}
Let $\phi:\Sigma\rightarrow M$ be an immersion from an oriented surface into an oriented four-dimensional Riemannian manifold $M$. Then the following assertions are equivalent:
\begin{description}
\item [i)] The twistor liftings $\tilde{\phi}^{\pm}:(\Sigma,J^{\pm})\rightarrow ({\cal P}^{\pm},{\cal J}^{\pm})$ of $\phi$ are holomorphic.
\item [ii)] The second fundamental form $\sigma$ of $\phi$ satisfies
\[
\sigma(u,v)=\sigma(J^{\pm}u,J^{\pm}v)-J^{\pm}\sigma(J^{\pm}u,v)-J^{\pm}\sigma(u,J^{\pm}v)
\]
for any vectors $u,v\in T\Sigma$.
\item [iii)] The almost complex structures $J^{\pm}$ define on $\phi^*TM$ structures of holomorphic bundles.
\end{description}
\end{proposition}

An immersion $\phi:\Sigma\rightarrow M$ satisfying one of the three equivalent conditions given in Proposition 2 will be called {\em twistor holomorphic with positive or negative spin}. As consequence of Proposition 2, the twistor holomorphicity of $\phi$ does not depend on the chosen orientation on the surface $\Sigma$. Hence we can talk about twistor holomorphic immersions from an {\em orientable} surface into an oriented four-dimensional Riemannian manifold. Twistor holomorphic surfaces with positive or negative spin which are also minimal are calles {\em superminimal surfaces with positive or negative spin}, ([B],[F],[G]). The surfaces which are simultaneously twistor holomorphic with positive and negative spin are the umbilical ones.

We define bilinear forms $\sigma^{\pm}$ on $\Sigma$ valuated on $T^{\perp}\Sigma$ by
\[
\sigma^{\pm}(u,v)=
\sigma(u,v)-\sigma(J^{\pm}u,J^{\pm}v)+J^{\pm}\sigma(J^{\pm}u,v)+J^{\pm}\sigma(u,J^{\pm}v)
\]
for vectors $u,v\in T\Sigma$. Then it is easy to check that
\begin{eqnarray*}
&{\displaystyle |H|^2+\bar K-\bar K^{\perp}=K-K^{\perp}+{1\over 16}|\sigma^+|^2,}& \\
&{\displaystyle |H|^2+\bar K+\bar K^{\perp}=K+K^{\perp}+{1\over 16}|\sigma^-|^2.}& 
\end{eqnarray*}
Now Proposition 2 gives the following known result ([F]), which relates the above functionals $W^{\pm}$ with the twistor theory.
\begin{proposition}
Let $\phi:\Sigma\rightarrow M$ be an immersion from an orientable compact surface into an oriented four-dimensional Riemanmnian manifold $M$. Then
\begin{description}
\item [i)] $W^+(\phi)\geq 2\pi(\chi-\chi^{\perp})$,

\item [ii)] $W^-(\phi)\geq 2\pi(\chi+\chi^{\perp})$.
\end{description} 
Moreover the equality in (i) holds if and only if $\phi$ is twistor holomorphic with positive spin, and the equality in (ii) holds if and only if $\phi$ is twistor holomorphic with negative spin.
\end{proposition}

In [ChT], $\chi-\chi^{\perp}$ was called the adjunction number of $\Sigma$ in $M$.

Perhaps the first case to study the functionals $W^{\pm}$ was when $(M,\langle ,\rangle)$ is the 4-dimensional Euclidean space, or equivalently (remember that $W^{\pm}$ are invariant under conformal transformation of the ambient space) when $(M,\langle ,\rangle )$ is a sphere $\S^4$ with its standard metric of constant curvature one. In this case, if $\phi:\Sigma\rightarrow \S^4$ is an immersion of an orientable compact surface, then it is easy to check that
\[
W^+(\phi)=W^-(\phi)=W(\phi)=\int_{\Sigma}\left(|H|^2+1\right)dA,
\]
and so these functionals are the classical Willmore functional $W$.

Moreover, as $\S^4$ is self-dual and anti-self-dual, the twistor spaces $({\cal P}^+,{\cal J}^+)$ and $({\cal P}^-,{\cal J}^-)$ are complex manifolds, and it is well-known that they are biholomorphic to $\C\P^3$ with its standard complex structure. Also the twistor projections $\pi^{\pm}$ are related by $\pi^-=A\circ\pi^+$, where $A$ is the antipodal map on $\S^4$, and hence the twistor holomorphic surfaces with negative spin of $\S^4$ are the images by the antipodal map of the twistor holomorphic surfaces with positive spin. In this case it is remarkable to refer to the papers [LY],[MoR] where lower bounds for $W$ are studied, and also to the paper [Mo] where the critical surfaces of genus zero of the functional $W$ are classified.

\section{Case of complex projective plane}
In $\C^{3}$ we consider the Hermitian product
\[
(z,w)=\sum_{i=1}^3z_i\bar w_i,
\]
for any $z,w\in\C^{3}$, where $\bar z$ stands for the conjugate of $z$. Then, $\Re (\,,\,)$ is the Euclidean metric and $\Im (\,,\,)$ the K\"{a}hler two--form on $\C^{3}$.

Let $\C\P^2$ be the complex projective plane with its canonical Fubini-Study metric $\langle ,\rangle $ of constant holomorphic sectional curvature $4$. Then
\[
\C\P^2=\{\Pi (z)=[z]\,/\, z=(z_1,z_2,z_3)\in\C ^{3}-\{0\}\}
\]
where $\Pi:\C^{3}-\{0\}\rightarrow\C\P^2$ is the standard projection. The metric $\Re (\,,\,)$ becomes $\Pi$ in a Riemannian submersion. The complex structure of $\C^{3}$ induces via $\Pi$ the canonical complex structure $J$ on $\C\P^2$. The K\"{a}hler two form $\Omega$ in $\C\P^2$ is defined by $\Omega(u,v)=\langle Ju,v\rangle $. We will consider $\C\P^2$ with the orientation $\Omega\wedge\Omega$.

If $\phi:\Sigma\rightarrow \C\P^2$ is an immersion of an oriented surface $\Sigma$, the {\em K\"{a}hler function} $C$ on $\Sigma$ is defined by
\[
\phi^*\Omega=C\,dA
\]
where $dA$ is the volume form on $\Sigma$.
We remark that the sign of $C$ depends on the orientation in $\Sigma$. So $C^2$ does not depend of the chosen orientation in $\Sigma$ and $C^2$ is defined even for non-orientable surfaces. It is clear that the K\"{a}hler function satisfies $-1\leq C\leq 1$. Surfaces with $C=1$, $C=-1$ and $C=0$ are called respectively {\em holomorphic}, {\em anti-holomorphic} and {\em Lagrangian}. A {\em complex surface} will be synonym of either a holomorphic or anti-holomorphic surface.

In addition, if $\Sigma$ is compact, then the {\em topological degree} $d$ of the map $\phi$ is given by
\[
d={1\over \pi}\int_{\Sigma}CdA.
\]
Also it is interesting to remark that the relation between $J$ and the almost complex structures $J^{\pm}$ on $\phi^*T\C\P^2$ defined in section 2 is
\begin{equation}
(Jv)^{\top}=CJ^+v =CJ^-v,\quad (J\xi)^{\perp}=CJ^+\xi=-CJ^-\xi,
\end{equation}
for any $v\in T\Sigma$ and $\xi\in T^{\perp}\Sigma$, where $\top$ and $\perp$ stand for tangent and normal components.

Before studying the functionals $W^{\pm}$ in this case, we are going to point out a strong property of the function $C$, which is not true when the codimension of the surface is bigger than two and that we have seen proved only in particular cases.
\begin{proposition}
Let $\Sigma$ be a compact orientable surface of $\C\P^2$ with constant K\"{a}hler function $C$. Then $\Sigma$ is either a complex or a Lagrangian surface.
\end{proposition}
{\em Proof:\/}
Suppose that $\Sigma$ is not a complex surface, i.e. the constant $C$ satisfied $-1<C<1$. Then, using (1), we can choose an oriented orthonormal local reference $\{e_1,e_2,e_3,e_4\}$ such that $\{e_1,e_2\}$ is an oriented reference on $T\Sigma$ and 
\begin{equation}
Je_1=Ce_2+\sqrt{1-C^2}e_4,\quad Je_2=-Ce_1+\sqrt{1-C^2}e_3.
\end{equation}
As $C=\langle Je_1,e_2\rangle$ and $C$ is constant, then derivating $C$ with respect to a tangent vector $v$ we have
\[
\langle \sigma(v,e_1),Je_2\rangle =\langle \sigma(v,e_2),Je_1\rangle .
\]
Using (2) in this formula we get
\[
A_{e_3}e_1=A_{e_4}e_2.
\]
Now it is straighforward to check, using the above information that $6C^2-K+K^{\perp}=0$. Integrating this equation, we finally obtain
\begin{equation}
6C^2\hbox{\,Area\,}(\Sigma)=2\pi(\chi -\chi^{\perp}).
\end{equation}
On the other hand, as $-1<C<1$, $F:T\Sigma\rightarrow T^{\perp}\Sigma$ defined by $F(v)=(Jv)^{\perp}$ defines an isomorphism of vector bundles, which implies that $\chi=\chi^{\perp}$. Using this in (3) we obtain that $C=0$, which finishes the proof.\hfill
q.e.d.

By using (1) and the well-known expression of the curvature of the Fubini--Study metric, it is straightforward to check that these functionals $W,W^{\pm}$ are given in this case by
\begin{eqnarray*}
 &{\displaystyle W(\phi)=\int_{\Sigma}\left(|H|^2+1+3C^2\right)dA}& \\
&{\displaystyle W^+(\phi)=\int_{\Sigma}\left(|H|^2+6C^2\right)dA, \quad W^-(\phi)=\int_{\Sigma}\left(|H|^2+2\right)dA.}&
\end{eqnarray*}
We remark that these functionals are defined for any compact surface not necessarily orientable.

As it was showed in section 2, to study the functionals $W^{\pm}$ it is necessary to understand the surfaces wich are twistor holomorphic. From the definition and using (1), it is easy to check that complex surfaces and minimal Lagrangian surfaces of $\C\P^2$ are twistor holomorphic surfaces with positive spin. In fact (as Gauduchon pointed out in [G]) they exactly are the superminimal surfaces with positive spin of $\C\P^2$. On the other hand, a complex or Lagrangian  twistor holomorphic surface with negative spin, must be an umbilical surface and then from [KZ] it must be totally geodesic. In the next result we study some important properties of these surfaces.
Before to stablish it, we need to point out a result which was proved in [EGT].

\begin{lemma} [EGT] Let $(\Sigma,\langle ,\rangle)$ be an oriented compact Riemannian surface, and $h:\Sigma\rightarrow\R$ a function of absolute value type, i.e., a smooth function satisfying $h=|t|f$ with $t$ a holomorphic function and $f$ a smooth positive function. Then
\[
\int_{\Sigma}\Delta\log h\, dA=-2\pi N(h),
\]
where $\Delta$ is the Laplacian operator of $\Sigma$ and $N(h)$ is the sum of all orders for all zeroes of $h$.
\end{lemma}

\begin{theorem}
Let $\phi:\Sigma\rightarrow\C\P^2$ be a twistor holomorphic immersion of an oriented surface $\Sigma$. If $\phi$ has positive spin, then
\begin{description}
\item [{i)}] Either $\phi$ is a complex immersion or the complex points of $\phi$ are isolated. Moreover the functions $\sqrt{1-C}$ and $\sqrt{1+C}$ are of absolute value type. 
\item [{ii)}] If $\Sigma$ is compact and non-complex, the degree $d$ of $\phi$ is given by
\[
3d=N_+-N_-,
\]
where $N_+=N(\sqrt{1-C})$ and $N_-=N(\sqrt{1+C})$.
\item [{iii)}] Under the conditions of $ii)$, 
\[
W^+(\phi)=\int_{\Sigma}\left(|H|^2+6C^2\right)dA=2\pi (N_++N_-).
\]
As consequence, if $\phi$ is totally real (i.e. $\phi$ has not complex points), then $\phi$ is a minimal Lagrangian surface.
\end{description}
If $\phi$ has negative spin, then
\begin{description}
\item [{i)}] The mean curvature $H$ is a holomorphic vector field on $T^{\perp}\Sigma$ with respect to the holomorphic structure associated to $J^-$. Hence either $\phi$ is superminimal with negative spin or $H$ has only isolated zeroes.
\item [{ii)}] If $\Sigma$ is compact and not superminimal, then
\[
W^-(\phi)=\int_{\Sigma}\left(|H|^2+2\right)dA=2\pi(\chi + N(H)),
\]
where $N(H)$ is the number of zeroes of $H$.
\end{description}
Moreover, if $\phi:\Sigma\rightarrow\C\P^2$ is an immersion of a sphere whose mean curvature $H$ is a holomorphic vector field on $T^{\perp}\Sigma$ with the holomorphic structure associated to $J^-$, then $\phi$ is either a complex immersion or a twistor holomorphic immersion with negative spin.
\end{theorem}

{\em Proof:\/} To prove the result we need to use complex coordinates on $\Sigma$. Let $z=x+iy$ be a local isothermal parameter on $\Sigma$ compatible with the given orientation. We will denote
\[
\partial\equiv \partial_z={1\over 2}(\partial_x-i\partial_y),\quad \bar{\partial}\equiv \partial_{\bar z}={1\over 2}(\partial_x+i\partial_y),
\]
the 	Cauchy-Riemann operators. Then
\[
|\partial_z|^2=\langle \partial_z,\partial_{\bar z}\rangle >0,\quad \langle \partial_z,\partial_ z\rangle =0,
\]
where $\langle,\rangle$  also denote the $\C$-linear extension of the metric $\langle ,\rangle $ to the complexified bundles. Then
\begin{equation}
\nabla_{\partial_z}\partial_z=\partial\log|\partial_z|^2\partial_z,\quad \sigma(\partial_z,\partial_{\bar z})=|\partial_z|^2H.
\end{equation}
If $\{e_3,e_4\}$ is an orthonormal local reference on $T^{\perp}\Sigma$ such that $\{\partial_x,\partial_y,e_3,e_4\}$ is the orientation on $\phi^*T\C\P^2$, we define
\[
\xi={e_3-ie_4\over \sqrt 2}.
\]
Then one can check (translating to complex notation the above arguments) that:
a) $\phi$ is twistor holomorphic with positive spin if and only if $\langle \sigma(\partial_z,\partial_z),\xi\rangle =0$;
b) $\phi$ is twistor holomorphic with negative spin if and only if $\langle \sigma(\partial_z,\partial_z),\bar{\xi}\rangle =0$; and
c) $J\partial_z=iC\partial_z+\langle J\partial_z,\xi\rangle \bar{\xi}$.

Suppose that $\phi$ is twistor holomorphic with positive spin. From a) we have that 
\begin{equation}
\sigma(\partial z,\partial z)=\langle\sigma(\partial z,\partial z),\bar{\xi}\rangle \xi.
\end{equation}
Then if $F=\langle J\partial_z,\xi\rangle $, using (4) and (5) it is clear that
\[
\partial F=\left(\partial\log|\partial_z|^2+\langle \nabla_{\partial_z}^{\perp}\xi,\bar{\xi}\rangle\right)F.
\]
So either $F$ vanishes identically or $F$ has only isolated zeros. But from c) we obtain that $|F|^2=|\partial_z|^2(1-C^2)$, and hence $\phi$ is either a complex immersion or the complex points of $\phi$ are isolated. This proves the first part of i).

On the other hand, the K\"{a}hler function $C$ is written by
\[
C={i\langle J\partial_{\bar z},\partial_z\rangle \over |\partial_z|^2}.
\]
So, using (4) and (5) again, it is easy to see that
\begin{equation}
\partial C=i\langle JH,\partial_z\rangle =-iF\langle H,\bar{\xi}\rangle .
\end{equation}
So derivating (6), using again (4) and (5) and the Codazzi equation, one can see that
\begin{equation}
\bar{\partial}\partial C=|\partial_z|^2(-|H|^2+3C(1-C^2)).
\end{equation}
Now, from (6) and (7), the gradient $\nabla C$ and the Laplacian $\Delta C$ of the fuction $C$ satisfied
\[
|\nabla C|^2=(1-C^2)|H|^2,\quad \Delta C=2C(-|H|^2+3(1-C^2)).
\]
If $\phi$ is not a complex immersion and from i) the complex points are isolated, then, outside the complex points, one can obtain easily from the above formulae
\[
\Delta\log \sqrt{1-C}=-{|H|^2\over 2}-3C(1+C),\quad \Delta\log \sqrt{1+C}=-{|H|^2\over 2}+3C(1-C).
\]
Now, ii) and iii) follow from Lemma 1 by proving that $\sqrt{1-C}$ and $\sqrt{1+C}$ are of absolute value type. This last assertion follows using a similar reasoning as in [EGT], Theorem A,(i), and its proof will be omited.

Suppose now that $\phi$ is twistor holomorphic with negative spin. From $b)$ we have that 
\begin{equation}
\sigma(\partial z,\partial z)=\langle(\partial z,\partial z),\xi\rangle \bar{\xi}.
\end{equation}
In complex notation, the holomorphicity of $H$ with respect to the holomorphic structure associated to $J^-$ means that
\begin{equation}
\nabla_{\partial_{\bar z}}^{\perp}(H-iJ^-H)=0.
\end{equation}
Derivating the second equation of (4), using (4), (8) and the Codazzi equation, it is not difficult to obtain that
\[
\nabla_{\partial_{\bar z}}^{\perp}H=B\xi,
\]
for a certain complex function $B$. From here it is clear that $H$ satisfied equation (9) and hence $H$ is a holomorphic vector field. The remaining assertions in $i)$ and $ii)$ are easy consequences of this fact.

Finally, let $\phi:\Sigma\rightarrow\C\P^2$ be an immersion of a sphere whose mean curvature vector $H$ is a holomorphic vector field with the holomorphic structure associated to $J^-$. We consider the very well-known cubic differential form $\Theta$ on $\Sigma$ (see [ES],[EGT]) defined by
\[
\Theta=\langle \sigma(\partial_z,\partial_z),J\partial_z\rangle dz^3.
\]
Now using similar arguments to used in the proof of the positive case, it is not difficult to check that the holomorphicity of $H$ implies that $\Theta$ is holomorphic, and then as $\Sigma$ is a sphere, $\Theta$ vanishes identically. So 
\begin{equation}
0=\langle \sigma(\partial_z,\partial_z),J\partial_z\rangle =\langle \sigma(\partial_z,\partial_z),\bar{\xi}\rangle F.
\end{equation}
As $H$ is holomorphic, either $H\equiv 0$ or $H$ has only isolated zeroes. If $H\equiv 0$, then (see [EGT]) either $\phi$ is a complex immersion or the complex points of $\phi$ are isolated. So from (10), if $\phi$ is not a complex immersion, $\phi$ is twistor holomorphic with positive spin and as consequence superminimal with positive spin. If $H$ has only isolated zeroes, the set of points where $F$ vanishes has empty interior, and so from (10) $\phi$ is twistor holomorphic with positive spin. This finishes the proof.\hfill q.e.d.

Now we are going to get the Euler-Lagrange equations for the functionals $W^+$ and $W^-$.
\begin{proposition}
Let $\phi:\Sigma\rightarrow\C\P^2$ be an immersion of an compact surface $\Sigma$. 
\begin{description}
\item [{i)}] $\phi$ is a critical point of the functional $W^+$ if and only if the mean curvature vector $H$ of $\phi$ satisfies
\[
\Delta^{\perp} H+(5+9C^2-2|H|^2)H+\tilde A (H)+12(JJ^+\nabla C)^{\perp}=0.
\]
\item [{ii)}] $\phi$ is a critical point of the functional $W^-$ if and only if the mean curvature vector $H$ of $\phi$ satisfies
\[
\Delta^{\perp} H+(1-3C^2-2|H|^2)H+\tilde A (H)=0.
\]
\end{description}
In both cases $\Delta^{\perp}$ and $\tilde A$ are defined by  
\[
\Delta^{\perp}=\sum_{i=1}^2\{\nabla^{\perp}_{e_i}\nabla^{\perp}_{e_i}-\nabla^{\perp}_{\nabla_{e_i}e_i}\},\quad 
\tilde A H=\sum_{e=1}^2\sigma(A_He_i,e_i),
\]
being $\{e_1,e_2\}$ an orthonormal reference tangent to $\Sigma$.
\end{proposition}

\begin{remark}
{\rm We note that minimal surfaces of $\C\P^2$ are critical points of the Willmore functional $W^-$.
However the only minimal surfaces critical for $W^+$ are the superminimal with positive spin, i.e., the complex and minimal Lagrangian surfaces [G]. In fact, from i) a minimal surface is critical for $W^+$ if and only if $JJ^+\nabla C$ is tangent to the surface. So $\Sigma_0=\{p\in \Sigma\,/\,(\nabla C)(p)\ne 0\}$ is an open subset of $\Sigma$ where $\phi$ is a complex immersion, which is imposible by the very definition of $\Sigma_0$. So $\Sigma_0=\o$, and then Proposition 4 proves the assertion.

On the other hand, notice that twistor holomorphic immersions with positive or negative spin are a kind of critical surfaces for $W^+$ or $W^-$ respectively because they minimize the corresponding Willmore functionals.}
\end{remark}

{\em Proof of Proposition 5:\/}
Following the computations got by Weiner in [W1], Theorem 2.1, it is easy to see that the first derivative of the functional $W^-$ is given by
\begin{eqnarray}
\delta W^-(\phi)=\int_{\Sigma}\langle \Delta H+\tilde{A}H+(1-3C^2-2|H|^2),\delta\phi\rangle\,dA,
\end{eqnarray}
where $\delta\phi$ stands for the variation vector field, which can be taken normal to the surface $\Sigma$. Now in order to compute the first derivative of the functional $W^+$, we start studying the functional $$
F(\phi)=\int_{\Sigma}C^2\,dA.$$ 

Using the well-known fact that
\begin{equation}
\delta(dA)=-2\langle H,\delta\phi\rangle\,dA,
\end{equation} 
it follows that the first derivative for the functional $F$ is given by
\[
\delta F(\phi)=\int_{\Sigma}\left(2C\delta(C)-2C^2\langle H,\delta\phi\rangle\right)\,dA.
\]
In order to compute $\delta(C)$, we recall the definition of $C$:
\[
\phi^*\Omega=C\,dA.
\]
Then, taking derivatives and using (12)
\begin{eqnarray}
\delta(\phi^*\Omega)=\delta(C)\,dA-2C\langle H,\delta\phi\rangle\, dA.
\end{eqnarray}
Now following standard arguments it is easy to check that if $\{e_1,e_2\}$ is an oriented orthonormal reference on $T\Sigma$, 
\[
(\delta(\phi^*\Omega))(e_1,e_2)=\langle \bar{\nabla}_{e_2}\delta\phi,e_1\rangle -\langle \bar{\nabla}_{e_1}\delta\phi,e_2\rangle =-\hbox{div\,} J^+(J\delta\phi)^{\top},
\]
where div stands for the divergence operator on $\Sigma$.
So, using this equation in (13) we get
\[
\delta(C)\,dA=\left(-\hbox{div\,} J^+(J\delta\phi)^{\top}+2C\langle H,\delta\phi\rangle\right)\,dA.
\]
So the first variation of $F(\phi)$ is 
\[
\delta F(\phi)=\int_{\Sigma}\left(2C^2\langle H,\delta\phi\rangle-2C\hbox{\,div\,}J^+(J\delta\phi)^{\top}\right)\,dA.
\]
On the other hand, from the divergence Theorem
\begin{eqnarray*}
&{\displaystyle
\int_{\Sigma}C\hbox{\,div\,}(J^+(J\delta\phi)^{\top}dA=-\int_{\Sigma}\langle \nabla C,J^+(J\delta\phi)^{\top}\rangle\,dA}&\\
&{\displaystyle
=-\int_{\Sigma}\langle (JJ^+\nabla C)^{\perp},\delta\phi\rangle\, dA.}&
\end{eqnarray*}
Using this we finally get
\[
\delta F(\phi)= \int_{\Sigma}2\langle C^2H+(J(J^+\nabla C)^{\perp},\delta\phi\rangle\,dA.
\]
Now using this formula, (11) and (12) we obtain the first variation of $W^+$, and the Proposition follows.\hfill q.e.d.

To finish this section, we are going to study twistor holomorphic surfaces from the view point of the twistor spaces.
If we consider $\C\P^2$ with the orientation $\Omega\wedge\Omega$, then $\C\P^2$ is a self-dual Riemannian manifold but not an anti-self-dual Riemannian manifold. So the twistor bundle $({\cal P}^-,{\cal J}^-)$ is a complex manifold and $({\cal P}^+,{\cal J}^+)$ is not a complex manifold. In fact it is well-known (see [ES]) that $({\cal P}^+,{\cal J}^+)$ can be differentiably identified with $P(T^{2,0}\C\P^2\oplus \C)$. Also, the complex manifold $({\cal P}^-,{\cal J}^-)$ can be endowed with a Riemannian metric which becomes it in a Einstein-K\"{a}hler manifold. Under this identification, ${\cal P}^-$ is the following complex hypersurface of $(\C\P^2\times \C\P^2,\langle,\rangle\oplus \langle ,\rangle ,J\oplus -J)$
\[
{\cal P}^-\equiv \{([z],[w])\in \C\P^2\times\C\P^2\,/\,z^t\bar w=0\},
\]
and the twistor projection $\pi^-$ is nothing but
\[
\pi^-(([z],[w]))=[\bar z\wedge \bar w],
\]
for any $[z],[w]\in \C\P^2$. Also the two natural projections $\pi_i:{\cal P}^-\rightarrow \C\P^2$ with $i=1,2$ are holomorphic and antiholomorphic maps respectively. 

The non-compact Lie group $PGL(3,\C)$ of the complex transformations of $\C\P^2$ acts over ${\cal P}^-$ by
\[
[A]\cdot ([z],[w])=([Az],[A^{*-1}w]),
\]
where $[A]$ is the class of a matrix $A\in GL(3,\C)$ and $A^*$ stands for the transpose conjugate of $A$. When $A\in U(3,\C)$, then
\[
\pi^-([A]\cdot ([z],[w]))=[A^t](\pi^-([z],[w])).
\]
\begin{remark}
{\rm Using this twistor space, one can reformulate, in a not very complicated way, the second part of Theorem 1 as follows.  {\em Given a non-complex immersion $\phi:\Sigma\rightarrow\C\P^2$ of an oriented surface $\Sigma$, then  $H$ is a holomorphic vector field in $T^{\perp}\Sigma$ with respect to the holomorphic structure associated to $J^-$ if and only if the twistor lifting $\tilde{\phi}$ of $\phi$ is a harmonic map}. So the second part of Theorem 1 is a consequence of the fact that every holomorphic map from an oriented surface into a K\"{a}hler manifold is harmonic.}
\end{remark}

The study (from this view point) of twistor holomorphic surfaces with positive spin is complicated because it is equivalent to study holomorphic curves in the non-complex manifold ${\cal P}^+$. Hovewer, as we will point out now, this view point will allow to understand very well the twistor holomorphic surfaces with negative spin. In fact,
if $\phi:\Sigma\rightarrow\C\P^2$ is a  twistor holomorphic immersion with negative spin, and $\tilde{\phi}=(\tilde{\phi}_1,\tilde{\phi}_2)$ its twistor lifting, then $\tilde{\phi}_1:\Sigma\rightarrow \C\P^2$ is a holomorphic curve and $\tilde{\phi}_2:\Sigma\rightarrow\C\P^2$ is an anti-holomorphic curve. In this context, $\phi$ is superminimal with negative spin if and only if $\tilde{\phi}_2$ is the dual curve of $\tilde{\phi}_1$. Hence, the Lie group $PGL(3,\C)$ acts on the twistor holomorphic surfaces with negative spin in the following way:
\[
[A]\cdot \phi=\pi^-([A]\cdot(\tilde{\phi}_1,\tilde{\phi}_2))=\pi^-([A]\tilde{\phi}_1,[A^{*-1}]\tilde{\phi}_2),
\]
for any $A\in GL(3,\C)$. This action sends superminimal surfaces with negative spin into themselves, and when $A\in U(3,\C)$, this action is the standard one. We will say that a twistor holomorphic surface with negative spin is a {\em twistor deformation} of another one if it is its image under the above action.

If $\Sigma$ is also compact, we will denote by $d_i$ the degree of $\tilde{\phi}_i$, $i=1,2$. If $d_i=0$ for some $i$, then $\tilde{\phi}_i$ is a point, and then (see definition of ${\cal P}^-$) $\tilde{\phi}_j$ with $j\not= i$ is a complex projective line. Hence $d_j=1$ and $\phi(\Sigma)$ is a complex projective line.

If $\phi$ is also a conformal immersion, then it is not difficult to check that $\tilde{\phi}_i$ are conformal immersions too and if $\langle,\rangle$  denotes the metric on $\Sigma$ induced by $\phi$ and $\langle ,\rangle _i$ the metrics on $\Sigma$ induced by $\tilde{\phi}_i$, $i=1,2$, then
\begin{equation}
\langle ,\rangle _1=\left({|H|^2+2(1-C)\over 4}\right)\langle ,\rangle ,\quad  \langle ,\rangle _2=\left({|H|^2+2(1+C)\over 4}\right)\langle ,\rangle .
\end{equation}
From here it follows the following result.

\begin{proposition}
Let $\phi:\Sigma\rightarrow \C\P^2$ be a twistor holomorphic immersion with negative spin of an orientable compact surface $\Sigma$. Then
\begin{description}
\item [{i)}] ${\displaystyle W^-(\phi)=\int_{\Sigma}\left(|H|^2+2\right)dA=2\pi( d_1+ d_2)}$,
\item [{ii)}] The degree $d$ of $\phi$ is given by $d=d_2 -d_1$.
\item [{iii)}] $W^-(\phi)\geq 2\pi$, and the equality holds if and only if $\phi(\Sigma)$ is a complex projective line $\C\P^1$.
\end{description}
\end{proposition}

In Proposition 6 above we have found that complex projective lines are twistor holomorphic suarfaces with negative spin attaining the minimum value for $W^-$. Now, we are going to describe other examples of twistor holomorphic surfaces with negative spin of $\C\P^2$ whose $W^-$ are also small.

{\em Twistor holomorphic compact surfaces with negative spin and $W^-=4\pi$}.

In this case, $(d_1,d_2)$ can take the value $(1,1)$. So, from Theorem 1 and Proposition 6, the surface must be a sphere with $\chi^{\perp}=0$ and $d=0$. From the remark made before Proposition 6, $\tilde{\phi}_2$ cannot be the dual curve of $\tilde{\phi}_1$, hence $\phi$ cannot be superminimal and then from Theorem 1, the holomorphic field $H$ has not zeroes. Also, as $\tilde{\phi}_i$ are unramified, (13) says that $\phi$ is unramified too. 

	Now, up to an holomorphic transformation of $\C\P^2$, $\tilde{\phi}_1:\C\cup\{\infty\}\rightarrow \C\P^2$ can be taken as
\[
\tilde{\phi}_1(z)=\Pi(1,z,0).
\]
Now, since $\tilde{\phi}=(\tilde{\phi}_1,\tilde{\phi}_2)$ lies in ${\cal P}^-$, easy computations say that $\tilde{\phi}_2:\C\cup\{\infty\}\rightarrow \C\P^2$ is given by
\[
\tilde{\phi}_2(z)=\Pi(\bar z,-1,\overline{P(z)}),\quad\hbox{with}\quad P(z)=a+bz\,,\,(a,b)\in\C^2-\{0\}.
\]
So our twistor holomorphic surface, (see definition of $\pi^-$), is a twistor deformation of $\phi_{a,b}:\C\cup\{\infty\}\rightarrow\C\P^2$, with $(a,b)\in\C^2-\{0\}$, where
\[
\phi_{a,b}(z)=\Pi(-a\bar{z}-b|z|^2,a+bz,1+|z|^2).
\]
It is interesting to remark that $\phi_{a,b}$ are embeddings.

{\em Twistor holomorphic compact surfaces with negative spin and $W^-=6\pi$.}

In this case, $(d_1,d_2)$ can take the value $(1,2),(2,1)$. But up to anti-holomorphic isometries of $\C\P^2$ it is sufficient to study the case $(d_1,d_2)=(1,2)$. In this case, the surface must be also a sphere but with $\chi^{\perp}=1$ and $d=1$. As $\tilde{\phi}_2$ cannot be the dual curve of $\tilde{\phi}_1$, $\phi$ cannot be superminimal and then from Theorem 1, the holomorphic field $H$ has only one zero. Also, as $\tilde{\phi}_i$ are unramified, (13) says that $\phi$ is unramified too. 

Now, as $PGL(3,\C)$ acts transitively on the conics of $\C\P^2$, $\tilde{\phi}_2:\C\cup\{\infty\}\rightarrow \C\P^2$ can be taken as
\[
\tilde{\phi}_2(z)=\Pi(1,\bar z,\bar z^2).
\]
So, since $\tilde{\phi}=(\tilde{\phi}_1,\tilde{\phi}_2)$ lies in ${\cal P}^-$, $\tilde{\phi}_1:\C\cup\{\infty\}\rightarrow \C\P^2$ is given by
\[
\tilde{\phi}_1(z)=\Pi(az,-a+bz,-b)\quad\hbox{with} \quad [(a,b)]\in\C\P^1.
\]
So our twistor holomorphic surface is a twistor deformation of $\psi_{[(a,b)]}:\C\cup\{\infty\}\rightarrow\C\P^2$, with $[(a,b)]\in\C\P^1$, where
\[
\psi_{[(a,b)]}(z)=\Pi(\bar b z(|z|^2+1)-\bar a z^2,-(\bar b+\bar a|z|^2z),\bar a(|z|^2+1)-\bar b\bar z).
\]

{\em Twistor holomorphic compact surfaces with negative spin and $W^-=8\pi$.}

The next examples are particulary interesting because they will be characterized in the next section. They will be twistor holomorphic spheres with negative spin, $W^-=8\pi$ and K\"{a}hler function $C=0$, i.e. Lagrangian surfaces. In particular their degrees will be zero, and then their twistor liftings will be a pair of conics. In fact these examples are called {\em Whitney spheres} and in [CU2] they were characterized as the only twistor holomorphic Lagrangian surfaces with negative spin. Up to isometries of $\C\P^2$ the are only a $1$-parameter family of surfaces which can be defined as
follow.

For each $t\in [0,\infty[$, we define $\phi_t:\S^2\rightarrow\C\P^2$ by
\[
\phi_{t}(x,y,z)=\Pi(x,y,z\cosh t +i\sinh t),
\]
for any $(x,y,z)\in\S^2=\{(x,y,z)\in\R^3\,/\,x^2+y^2+z^2=1\}$.

We remark that $\phi_0$ is the totally geodesic immersion of $\S^2$, which is a covering of the totally geodesic embedding of $\R\P^2$. For $t>0$, $\phi_t$ is an embedding except at the poles of $\S^2$ where it has a double point. Amongst them only $\phi_0$ is a minimal surface.

\section{Lower bounds for the functional $W^-$}
In this section we start obtaining a lower bound for the functional $W^-$.

\begin{theorem}
Let $\phi:\Sigma\rightarrow \C\P^2$ be an immersion of a compact surface $\Sigma$ and $\mu$ the maximum multiplicity of $\phi$. Then:
\[
\int_{\Sigma}\left(|H|^2+3+C^2\right)\,dA\geq 4\pi\mu,
\]
and the equality holds if and only if $\phi(\Sigma)$ is a complex projective line $\C\P^1$.
As a consequence,
\[
W^-(\phi)=\int_{\Sigma}\left(|H|^2+2\right)dA\geq 2\pi\mu,
\]
and the equality holds if and only if $\phi(\Sigma)$ is a complex projective line $\C\P^1$.
\end{theorem}

\begin{corollary}
The area of a compact minimal surface $\Sigma$ immersed in $\C\P^2$ with maximum multiplicity $\mu$ satisfies
\[
\hbox{{\rm Area}\,}(\Sigma)\geq \pi\mu,
\]
and the equality holds if and only if $\Sigma$ is a complex projective line of $\C\P^2$.
\end{corollary}

{\em Proof of Theorem 2:\/}
As the maximum multiplicity of $\phi$ is $\mu$, let $\{p_1,\dots,p_{\mu}\}$ be points of $\Sigma$ such that $\phi(p_i)=[a]\in\C\P^2$ for any $i=1,\dots,\mu$. We define a function $f:\C\P^2\rightarrow\R$ by
\[
f([z])={|(z,a)|^2\over |z|^2|a|^2},
\]
for any $[z]\in\C\P^2$. Then $0\leq f\leq 1$, and $f([z])=0$ if and only if $[z]$ is in the cut locus $\C\P^1_{[a]}$ of the point $[a]$. Also, $f([z])=1$ if and only if $[z]=[a]$. So $\log (1-f)$ is a well defined function on $\C\P^2-\{[a]\}$. 

From now on (in order to simplify the notation) we will consider that $|a|=1$, and we will restrict $\Pi$ to the unit sphere $\S^{5}\subset\C^3$. So, the function $f$ will be nothing but
\[
f([z])=|(z,a)|^2,
\]
for any $z\in\S^{5}$. 

First we compute the gradient of $f$. If $v$ is any tangent vector to $\C\P^2$ at $[z]$, then
\[
v(f)=2\Re (v^*,(z,a)a),
\]
being $v^*$ the horizontal lift to $T_z\S^5$ of $v$.
So, 
\[
(\bar{\nabla} f)_{[z]}=2(d\Pi)_z((z,a)a-|(z,a)|^2z),
\]
for any $[z]\in \C\P^2$. It is intereting to remark that $|\bar{\nabla} f|^2=4f(1-f)$.

Now, using that $\Pi:\S^{5}\rightarrow\C\P^2$ is a Riemannian submersion, it is easy to check that the Hessian of $f$ is given by
\begin{eqnarray}
(\bar{\nabla}^2f)(u,v)=-2f\langle u,v\rangle +2\Re ((u^*,a)(a,v^*)),
\end{eqnarray}
for any vectors $u,v\in T_{[z]}\C\P^2$, being $u^*,v^*$ the horizontal lifts to $T_z\S^5$ of $u,v$. In that follows it will be interesting to take into account the following formula, which can be easily check
\begin{equation}
f\Re ((u^*,a)(a,v^*))=\langle\bar{\nabla}f,u\rangle\langle\bar{\nabla}f,v\rangle+ \langle\bar{\nabla}f,Ju\rangle\langle\bar{\nabla}f,Jv\rangle
\end{equation}

We can define on $\Sigma-\{p_1,\dots,p_{\mu}\}$ the function $\log (1-h)$, where $h=f(\phi)$.

By decomposing $$\bar{\nabla}f\circ \phi=\nabla h+\xi$$ in its tangencial and normal components, it is easy to see that
\begin{eqnarray}
&|\nabla h|^2=4h(1-h)-|\xi|^2,&\nonumber\\& &  \\&(\nabla^2 h)(u,v)=(\bar{\nabla}^2f)(\phi_*u,\phi_*v)+\langle \sigma(u,v),\xi\rangle,&\nonumber
\end{eqnarray}
for any vectors $u,v$ tangent to $\Sigma$. From (17) and (15) we obtain that
\[
\Delta \log (1-h)=-{2\over 1-h}\sum_{i=1}^2|(e_i^*,a)|^2 
-{2\over 1-h}\langle H,\xi\rangle  +{|\xi|^2\over (1-h)^2},
\]
where $\{e_1,e_2\}$ is an orthonormal reference on $\Sigma$. As
\[
\left|H-{\xi\over 1-h}\right|^2=|H|^2+{|\xi|^2\over (1-h)^2}-{2\over 1-h}\langle H,\xi\rangle ,
\]
we obtain that
\begin{equation}
\Delta \log (1-h)\geq -|H|^2-{2\over 1-h}\sum_{i=1}^2|(e_i^*,a)|^2,
\end{equation}
and the equality holds if and only if $H=\xi/(1-h)$.

Now, from (16), we obtain
\begin{equation}
4h\sum_{i=1}^2|(e_i^*,a)|^2=(1+C^2)|\nabla h|^2+(1-C^2)|\xi|^2
+2\langle (J\xi)^{\perp},J\nabla h\rangle.
\end{equation}
On the other hand
\begin{eqnarray*}
& 0\geq -|\sqrt{2}(J\xi)^{\perp}-{1\over \sqrt{2}}(J\nabla h)^{\perp}|^2&\\
&=-2C^2|\xi|^2-{1\over 2}(1-C^2)|\nabla h|^2+2\langle (J\xi)^{\perp},J\nabla h\rangle.&
\end{eqnarray*}
Using this inequality in (19) we get
\begin{eqnarray*}
&{\displaystyle 4h\sum_{i=1}^2|(e_i^*,a)|^2\leq {3+C^2\over 2}|\nabla h|^2+(1+C^2)|\xi|^2}&\\
&{\displaystyle ={3+C^2\over 2}|\bar\nabla f\circ\phi|^2-{1-C^2\over 2}|\xi|^2\leq {3+C^2\over 2}4h(1-h),}&
\end{eqnarray*}
and the equality holds if and only if $(1-C^2)\xi=0$ and $(J\xi)^{\perp}$ and $(J\nabla h)^{\perp}$ are colinear.

So finally, from (18) we get that
\begin{equation}
\Delta \log (1-h)\geq -|H|^2-3-C^2,
\end{equation}
and the equality holds if and only $\phi(\Sigma)$ is a complex projective line $\C\P^1$. In fact, the equality holds if and only if $H=\xi/ (1-h)$, $(1-C^2)\xi=0$ and $(J\xi)^{\perp}$ and $(J\nabla h)^{\perp}$ are colinear. Let
\[
\Sigma_0=\{p\in\Sigma\,/\,C^2(p)=1\}.
\]
If the interior of $\Sigma_0$ is empty, then $\xi=0$ on the whole $\Sigma$. If not the interior of $\Sigma_0$ is a complex curve and then it is a minimal surface. So, as $\xi=(1-h)H$, we have $\xi=0$ in the interior of $\Sigma_0$. So in any case $\xi\equiv 0$ on the whole surface and in particular the surface is minimal. The third condition says that $J\nabla h$ is a tangent vector, and so outside the zeroes of $\nabla h$, $\phi$ is  a complex curve. On this set, and because $\bar{\nabla}f\circ\phi=\nabla h$, (16) says that $\sigma(v,\nabla h)=0$ for any vector $v$ tangent to $\Sigma$, and  so $\phi$ is totally geodesic. As $|\nabla h|^2=4h(1-h)$, 
\[
\{p\in\Sigma\,/\,(\nabla h)(p)=0\}=\{p_1,\dots,p_{\mu}\}\cup\phi^{-1}\left(\C\P^1_{[a]}\right).
\]
So $\phi$ is a complex totally geodesic surface on the whole $\Sigma$ which finishes our claim.

Let $B_{[a]}(\varepsilon)$ be the geodesic ball in $\C\P^2$ centered at the point $[a]$ with radius $\arccos\sqrt{1-\varepsilon^2}$, that is, the set$$
B_{[a]}(\varepsilon)=\{p\in\C\P^2\,|\,1-f(p)\le\varepsilon^2\},$$
with $\varepsilon$ too small in order to $B_\varepsilon=\phi^{-1}(B_{[a]}(\varepsilon))$ will be the disjoint union of neighbourhoods $B_i$, $i=1,\dots,\mu$  around $p_i$ in $\Sigma$. Then the divergence theorem on the manifold $\Sigma-B$  says that
\[
\int_{\Sigma-B}\Delta \log (1-h)\, dA=-\sum_{i=1}^{\mu}\int_{\partial B_i}\frac{\langle \nabla h,\nu_i\rangle}{1-h}\, ds,
\]
where $\nu_i$ is the unit conormal of $\partial B_i$ pointing to the interior of $B_i$. Since the function $h$ attains its maximum value 1 at each $p_i$ and $h$ is constant along each $\partial B_i$, we have that$$
\nu_i=\frac{\nabla h}{|\nabla h|}_{|\partial B_i}.$$
So, combining these equalities with the integral equality above, we have$$
\int_{\Sigma-B}\Delta \log (1-h)\, dA=-\sum_{i=1}^{\mu}\int_{\partial B_i}\frac{|\nabla h|}{1-h}\, ds=-\sum_{i=1}^{\mu}\frac{1}{\varepsilon^2}\int_{\partial B_i}|\nabla h|\, ds.$$
As $\varepsilon$ tends to zero, $|\nabla h|$ along $\partial B_i$ approaches to $|\nabla f|=
2\varepsilon\sqrt{1-\varepsilon^2}$ and the length of $\partial B_i$ approaches to$$
2\pi\hbox{\,\rm radius\,}B_i=2\pi\arccos\sqrt{1-\varepsilon^2}.$$
Then, we obtain that$$
\int_\Sigma\Delta\log (1-h)\,dA=-4\pi\mu.$$
This equality and (20) prove the inequality we were looking for. \hfill q.e.d.

In the next result we improve the lower bound obtained in Theorem 2 for $W^-$ in the family of Lagrangian surfaces of $\C\P^2$. 

\begin{theorem}
Let $\phi:\Sigma\rightarrow \C\P^2$ be a Lagrangian  immersion of a compact surface $\Sigma$ and 
$\mu$ the maximum multiplicity of $\phi$. Then
\[
W^-(\phi)=\int_{\Sigma}\left(|H|^2+2\right)\,dA\geq 4\pi\mu,
\]
and the equality holds if and only if either $\phi$ is totally geodesic and $\phi(\Sigma)$ is a  real projective plane with $W^-(\phi)=4\pi$ or $\phi$ is a Whitney sphere with $W^-(\phi)=8\pi$.
\end{theorem}

\begin{corollary}
The area of a minimal Lagrangian compact surface $\Sigma$ immersed in $\C\P^2$ with maximum multiplicity $\mu$ satisfies
\[
\hbox{{\rm Area}\,}
(\Sigma)\geq 2\pi\mu,
\]
and the equality holds if and only if $\Sigma$ is totally geodesic.
\end{corollary}

\begin{remark}
{\rm  We consider the following holomorphic an anti-holomorphic immersions $\tilde{\phi}_i:\C\cup\{\infty\}\rightarrow\C\P^2$, $i=1,2$ given by
\[
\tilde{\phi}_1(z)=\Pi(1,z,0),\quad \tilde{\phi}_2(z)=\Pi(\bar z(1+\bar z),-(1+\bar z),\bar z).
\]
Then the corresponding twistor holomorphic surface with negative spin $\phi:\C\cup\{\infty\}\rightarrow\C\P^2$ is 
\[
\phi(z)=\Pi(-|z|^2,z,(1+z)(1+|z|^2)).
\]
It is easy to check that $\phi$ is regular and that is embedded except at $z=0,\infty$ where $\phi$ has a doble point. So $\mu=2$. As the degrees of $\tilde{\phi}_1$ and $\tilde{\phi}_2$ are $1$ and $2$, then $W^-(\phi)=6\pi$. This example shows that even in the family of non-complex compact surfaces of $\C\P^2$, Theorem 3 is not true.}
\end{remark}

{\em Proof of Theorem 3:\/}
Using a similar reasoning like in the proof of Theorem 1, we get the following integral formula
\[
\int_{\Sigma}\left(|H|^2+{2\over 1-h}\sum_{i=1}^2|(e_i^*,a)|^2\right)\,dA\geq4\pi\mu,
\]
and the equality holds if and only if $H=\xi/(1-h)$. In this case, using that $\phi$ is a Lagrangian immersion, i.e., $C=0$, in (19) we have
\[
\sum_{i=1}^2|(e_i^*,a)|^2=1-h.
\]
So we finally obtain
\[
\int_{\Sigma}\left(|H|^2+2\right)\,dA\geq4\pi\mu,
\]
and the equality holds if and only if $H=\xi/(1-h)$.
Now we are going to classify Lagrangian surfaces of $\C\P^2$ whose mean curvature is given in the above way. 

From now on we will work on the dense open subset of $\Sigma$ defined by
\[
\Sigma_0=\{p\in\Sigma\,/\,\bar\nabla f(\phi(p)) \neq 0.\}
\]
On $\Sigma_0$ the function $h$ satisfies $0<h<1$. First, from (17), (15) and using elementary properties of Lagrangian surfaces it follows
\begin{equation}
\langle\nabla_vJ\xi,w\rangle=\langle\sigma(v,w),J\nabla h\rangle-2\Re ((v^*,a)(a,(Jw)^*)).
\end{equation}
So derivating $JH=J\xi/(1-h)$ and using (16) and (21) it follows
\begin{eqnarray*}
&{\displaystyle \langle \nabla_v JH,w\rangle ={1\over 1-h}\langle \sigma(v,w),J\nabla h\rangle +{2-h\over h(1-h)^2}\langle \nabla h,v\rangle \langle J\xi,w\rangle }&\\
&{\displaystyle -{2\over h(1-h)}\langle \nabla h,w\rangle \langle J\xi,v\rangle,}&
\end{eqnarray*}
for any $v,w$ tangent to $\Sigma$. As $JH$ is a closed vector field on $\Sigma$, the first term is symmetric. Also, as the second fundamental form is also symmetric, we obtain that the others terms are symmetric too, and so
\begin{eqnarray}
dh\wedge \alpha=0,
\end{eqnarray}
where $\alpha$ is the 1-form on $\Sigma$ given by
\[
\alpha(v\rangle =\langle v,J\xi\rangle .
\]
Now we are going to prove that there exists a vector field $X$ tangent to $\Sigma_0$ and functions $a$ and $b$ on $\Sigma_0$ with $a^2+b^2=h$ such that
\begin{equation}
\bar{\nabla}f\circ\phi=aX+bJX.
\end{equation}
So in particular this vector field $X$ verifies $|X|^2=4(1-h)$.

In fact, let $A=\{p\in\Sigma_0\,/\,dh_p=0\}$ and $B=\{p\in\Sigma_0\,/\,\alpha _p=0\}$. If $A=\Sigma_0$ or $B=\Sigma_0$ the claim is trivial. Otherwise, $A$ and $B$ are proper closed subsets of $\Sigma_0$. Now, $dh\wedge\alpha=0$ says that on $\Sigma_0/A$ we can writte $\alpha=\lambda dh$ for certain smooth function $\lambda$. Taking $X={\sqrt{1+\lambda^2}\over\sqrt h}\nabla h$, then on $\Sigma_0/A$ we have $\bar\nabla f\circ\phi=aX+bJX$ for certain smooth functions $a$ and $b$ on $\Sigma_0/A$ satisfying $a^2+b^2=h$. Making a similar reasoning with $B$ we writte on $\Sigma_0/B$, $\bar\nabla f\circ\phi=a'X'+b'JX'$ with $a'^2+b'^2=h$. It is clear that, on the non-empty subset $\Sigma_0/(A\cup B)$, we can take $X'=X$, $a'=a$ and $b'=b$. So we prove the existence of such $X$ satisfying (23).

Derivating (23) with respect to a vector $v$, taking tangent and normal components and using (15), (16) and (17) we obtain
\begin{eqnarray*}
&-2hv+2\langle v,X\rangle X=\langle \nabla a,v\rangle X+a\nabla_vX-bA_{JX}v& \\
&\langle \nabla b,v\rangle X+aA_{JX}v+b\nabla_vX=0.&
\end{eqnarray*}
From these equations it is easy to obtain that
\[
\sigma(X,X)=2\rho JX,\quad \sigma(X,V)=2b JV,
\]
for certain function $\rho$, where $V$ is any orthogonal vector field to $X$. From here
\[
2H={\rho +b\over 2(1-h)}JX +cJV,
\]
for certain function $c$. But $\xi=bJX$, and then 
\[
H={b\over 1-h}JX.
\]
So we get that $\rho=3b$ and $c=0$. This in particular means that $|\sigma|^2=3|H|^2$. The Gauss equation implies that $|H|^2+2=2K$, and since our surface is Lagrangian, $K^{\perp}=K$. So finally we get that $|H|^2+2=K+K^{\perp}$, which means that our surface is twistor holomorphic with negative spin. Now the mean result in [CU2] finishes our proof.\hfill q.e.d.

\begin{remark}
{\rm The totally geodesic surfaces and the Whitney spheres of $\C\P^2$ have the property that their mean curvature vectors are given by $$
H={(\bar{\nabla}f\circ\phi)^{\perp}\over 1-h},$$
being $h=f\circ\phi$, $f([z])=|(z,a)|^2$ and $\phi$ the immersion. The geometric meaning of this property is the following. Let $M=\C\P^2-{[a]}$, and $g$ the metric on $M$ conformal to the Fubini-Study metric defined by
\[
g={1\over (1-f)^2} \langle,\rangle.
\]
Then it is not difficult to see that $(M,g)$ is a complete Riemannian manifold (with one end) and with zero scalar curvature. If $\phi:\Sigma\rightarrow\C\P^2$ is an immersion with $\{p_1,\dots,p_{\mu}\}=\phi^{-1}([a])$, then the mean curvatures vectors $\hat{H}$ and $H$ of $\phi$ with respect to the metric induced by $g$ and $\langle,\rangle$ are related (see proof of Proposition 1) by
\[
{\hat{H}\over 1-f^2}=H-{(\bar{\nabla}f\circ\phi)^{\perp}\over 1-h}.
\]
So the condition $H=(\bar{\nabla}f\circ\phi)^{\perp} /(1-h)$ means that the surface $\Sigma-\{p_1,\dots,p_{\mu}\}$ is minimal in $(M,g)$.}
\end{remark}

To end this section we would like to remark something about the functional $W^-$ acting on tori.
When you consider the Willmore functional on compact surfaces of $\R^4$, there is a very famous conjecture, due to Willmore, which says that the Willmore functional on tori is bounded below by $2\pi^2$ and the Clifford torus is the only torus which achieves this minimum. Since the Clifford torus is Lagrangian, Minicozzi [M], studied this problem in this smaller class of Lagrangian tori.

In our case, we will also called {\em Clifford torus} to the following torus $T$ embedded in $\C\P^2$ and defined by
\[
T=\{\Pi(z)\in \C\P^2\,/\,|z_i|^2={1\over 3},i=1,2,3\}.
\]
It is easy to check that $T$ is a minimal Lagrangian torus with area $4\pi^2/ 3\sqrt 3$. So its Willmore functional is $W^-(T)=8\pi^2/ 3\sqrt 3$. 

For complex tori of $\C\P^2$, the Willmore functional $W^-$ take the value $2\pi d$, where $d$ is the degree of the torus. So, as $d$ must be non smaller that $3$, we obtain that in this family $W^-\geq 6\pi>8\pi^2/ 3\sqrt 3$. 

For twistor holomorphic tori in $\C\P^2$ with negative spin, Proposition 4 says that the Willmore functional $W^-$ satisfied again $W^-=2\pi(d_1+d_2)$. But again the corresponding holomorphic and antiholomorphic curves have genus one. So $W^-\geq 6\pi$.

In the family of tori of $\C\P^2$ with non-zero parallel mean curvature vector, very recently Kenmotsu and Zhou in [KZ] have proved that they are Lagrangian and flat, and then, up to isometries they can be parametrized by $T_{r_1,r_2,r_3}\,/\,r_1\geq r_2\geq r_3>0,r_1^2+r_2^2+r_3^2=1$, where
\[
T_{r_1,r_2,r_3}=\{\Pi(z)\in\C\P^2\,/\;|z_i|^2=r_i^2\,,\, i=1,2,3\}.
\]
It is clear that $W^-(T_{r_1,r_2,r_3})\geq{8\pi^2 r_1r_2r_3}\geq 8\pi^2/ 3\sqrt 3$, and the equality is only achieved by the Clifford torus.

Also, in [CU1], Castro and Urbano classified minimal Lagrangian tori of $\C\P^2$ invariant under a $1$-parameter group of holomorphic isometries. This family of tori is described in terms of elliptic functions, and it is not a complicated exercise to check that the Willmore functional $W^-$ on these tori satisfy $W^-\geq 8\pi^2/ 3\sqrt 3$, with equality only for the Clifford torus.

These considerations make reasonable the following conjecture:
\begin{quote}
{\em The Clifford torus achieves the minimum of the Willmore functional $W^-$ either amongst all tori in $\C\P^2$ or amongst all Lagrangian tori in $\C\P^2$}.
\end{quote}

\begin{tabular}{ll} 
Departamento de Geometr\'{\i}a y Topolog\'{\i}a \\
Universidad de Granada \\
18071 Granada \\
SPAIN \\
e-mails:smontiel@goliat.ugr.es  and  furbano@goliat.ugr.es \\
\end{tabular}
\end{document}